\documentclass[10pt]{amsart}
\usepackage{amsmath,amsthm,amssymb,amsfonts,amscd}
\usepackage{epsfig}
\usepackage{xypic}
\numberwithin{equation}{section}

\theoremstyle{plain}
\newtheorem{theorem}{Theorem}

\newtheorem*{theorem*}{Theorem}
\newtheorem*{conjecture*}{Conjecture}

\theoremstyle{definition}
\newtheorem{remark}[theorem]{Remark}

\newcommand{\CC}{{\mathbb{C}}}
\newcommand{\HH}{{\mathbb{H}}}

\newcommand{\FF}{{\mathbb{F}}}
\newcommand{\RR}{{\mathbb{R}}}
\newcommand{\ZZ}{{\mathbb{Z}}}

\begin{document}
\title[A McKay correspondence for some finite subgroups of ${\rm SL}_3(\CC)$]{A McKay correspondence for the Poincar\'e series of some finite subgroups of ${\rm SL}_3(\CC)$}
\author{Wolfgang Ebeling}
\thanks{Partially supported by DFG}
\address{Institut f\"ur Algebraische Geometrie, Leibniz Universit\"at Hannover, Postfach 6009, D-30060 Hannover, Germany}
\email{ebeling@math.uni-hannover.de}
\dedicatory{Dedicated to the memory of Egbert Brieskorn with great admiration}
\date{}
\begin{abstract} A finite subgroup of ${\rm SL}_2(\CC)$ defines a (Kleinian) rational surface singularity. The McKay correspondence yields a relation between the Poincar\'e series of the algebra of invariants of such a group and the characteristic polynomials of certain Coxeter elements determined by the corresponding singularity. Here we consider some non-abelian finite subgroups $G$ of ${\rm SL}_3(\CC)$. They define non-isolated three-dimensional Gorenstein quotient singularities. We consider suitable hyperplane sections of such singularities which are Kleinian or Fuchsian surface singularities. We show that we obtain a similar relation between the group $G$ and the corresponding surface singularity.  
\end{abstract}
\maketitle
\section*{Introduction}

In \cite{Ebeling02} we showed that the Poincar\'e series of the coordinate algebra of a two-dimensional quasihomogeneous singularity is the quotient of two polynomials one of which is related to the characteristic polynomial of the monodromy of the singularity. There are two special cases of this result. One is the case of a Kleinian singularity not of type $A_{2n}$. The Kleinian singularities are defined by quotients of $\CC^2$ by finite subgroups of ${\rm SL}_2(\CC)$. In this case, the relation means that the Poincar\'e series is the quotient of the characteristic polynomials of the Coxeter element and the affine Coxeter element of the corresponding root system of type ADE. We derived this relation from the McKay correspondence. The other case is the case of a Fuchsian singularity. A Fuchsian singularity is defined by the action of a Fuchsian group (of the first kind) $\Gamma \subset {\rm PSL}(2, \RR)$ on the tangent bundle $T_\HH$ of the upper half plane $\HH$. For a Fuchsian hypersurface singularity (or more generally for a Fuchsian singularity of genus 0 \cite{EP1}), we showed that the Poincar\'e series is the quotient of two characteristic polynomials of Coxeter elements \cite{Ebeling03}. 

Here we consider a similar relation for the Poincar\'e series of some non-abelian finite subgroups of ${\rm SL}_3(\CC)$. The non-abelian finite subgroups of ${\rm SL}_3(\CC)$ define non-isolated three-dimensional Gorenstein quotient singularities. We consider those groups where the natural three-dimensional representation is irreducible and the corresponding quotient singularity has a certain hyperplane section which is a Kleinian or Fuchsian singularity.  We show, that in this way, we again obtain relations between the Poincar\'e series of the algebra of invariants of the group and the characteristic polynomials of certain Coxeter elements determined by the corresponding Kleinian or Fuchsian singularity. 

The famous paper \cite{Br} of E.~Brieskorn is fundamental for the study of Kleinian singularities. The Kleinian singularities were a central theme in Brieskorn's research and we owe Brieskorn many beautiful and important results about these singularities. Therefore I would like to express my great admiration for him in dedicating this paper to his memory. 

\section{Finite subgroups of ${\rm SL}_2(\CC)$ and ${\rm SL}_3(\CC)$ and normal surface singularities}
Let $G$ be a finite subgroup of ${\rm SL}_2(\CC)$. The classification of finite subgroups of ${\rm SL}_2(\CC)$ up to linear equivalence is well-known, see e.g.\ \cite{Kle93}. There are up to conjugacy five classes of such groups: the cyclic groups $\mathcal{C}_{\ell}$, the binary dihedral groups $\mathcal{D}_n$, the binary tetrahedral group $\mathcal{T}$, the binary octahedral group $\mathcal{O}$, and the binary icosahedral group $\mathcal{I}$. The quotients of $\CC^2$ by these groups were studied by E.~Brieskorn \cite{Br}. Equations for these singularities can be obtained from generators and relations of the algebra of invariant polynomials with respect to $G$. This algebra has three generators $x,y,z$  in each case which satisfy an equation $R(x,y,z)=0$. The degrees of the generators and the equation $R(x,y,z)=0$ are indicated in Table~\ref{TabSL2sing}. (They can be found, e.g., in \cite{Sp}.) The equations define isolated hypersurface singularities in $\CC^3$, the so called {\em Kleinian singularities}. 
\begin{table}
\begin{center}
\begin{tabular}{|c|c|c|c|c|c|c|}
\hline
 $G$    & $|G|$ & $x,y,z$ & $c_G$ & $R(x,y,z)$ & Sing. & $\alpha_1, \ldots , \alpha_m$ \\
\hline
$\mathcal{C}_{2n+1}$ & $2n+1$ & $2,2n+1,2n+1$ & 1 & $x^{2n+1}+y^2+z^2$ & $A_{2n}$ & $2n$\\
$\mathcal{C}_{2n}$ & $2n$ & $2,2n,2n$ & 2 & $x^{2n}+y^2+z^2$ & $A_{2n-1}$ & $2n-1$\\
$\mathcal{D}_n$ & $4n$ & $4,2n,2n+2$ & 2 & $x^{n+1} +xy^2+z^2$ & $D_{n+2}$ & $2,2,n$\\
$\mathcal{T}$ & $24$ & $6,8,12$ & 2 & $x^4+y^3+z^2$ & $E_6$ & $2,3,3$ \\
$\mathcal{O}$ & $48$ & $8,12,18$ & 2 & $x^3y+y^3+z^2$ & $E_7$ & $2,3,4$ \\
$\mathcal{I}$ & $120$ & $12,20,30$ & 2 & $x^5+y^3+z^2$ & $E_8$ & $2,3,5$\\
\hline
\end{tabular}
\end{center}
\caption{Subgroups of ${\rm SL}_2(\CC)$ and surface singularities} \label{TabSL2sing}
\end{table}

\begin{sloppypar}

The finite subgroups of ${\rm SL}_3(\CC)$ were classified up to linear equivalence by H.~F.~Blichfeldt, G.~A.~Miller, and L.~E.~Dickson \cite{Bli,MBD} with two missing cases (see \cite{YY}). There are 12 types of finite subgroups of ${\rm SL}_3(\CC)$: (A)--(L). There are four infinite series (A)--(D). The groups of type (A) are the diagonal abelian groups and the groups of type (B) are isomorphic to transitive finite subgroups of ${\rm GL}_2(\CC)$.  Here the natural 3-dimensional representation is not irreducible. Type (C) is the infinite series $\Delta(3n^2)$ of groups and  type (D) the series $\Delta(6n^2)$ (for the notation see \cite{HH, LNR, EL}). Moreover, we have 8 exceptional subgroups (E)--(L). 

\end{sloppypar}
We consider those subgroups of type (C)--(L) which admit a certain hyperplane section which defines a Kleinian or Fuchsian singularity. 
Generators and relations for the algebra of invariant polynomials with respect to $G$ have been computed in \cite{YY} (see also \cite{W} for some cases). They correspond to non-isolated Gorenstein quotient singularities $\CC^3/G$. These singularities are either hypersurface singularities in $\CC^4$ or complete intersection singularities in $\CC^5$. We denote the coordinates of these spaces by $w,x,y,z$ and $w,x,y,z,u$ respectively. We consider hyperplane sections of these singularities, namely we consider the restrictions of the equations to the hyperplane $w=0$. For the series (C) and (D) the hyperplane sections of the corresponding singularities for the first few elements of these series are listed in Table~\ref{TabCD}. It turns out that the singularities corresponding to the series (C) ($\Delta(3n^2)$) belong to Arnold's $E$-series whereas those of type (D) ($\Delta(6n^2)$) belong to Arnold's $Z$-series ($n$ even) or are complete intersection singularities ($n$ odd) (for the definition of these series see \cite{Ar}). The subgroups which correspond to Kleinian singularities are the tetrahedral group $T=\Delta(3 \cdot 2^2)$ and the octahedral group $O=\Delta(6\cdot 2^2)$ which correspond to the Kleinian singularities $E_6$ and $E_7$ respectively. Those which correspond to Fuchsian singularities are $\Delta(3 \cdot 4^2)$ ($E_{14}$), $\Delta(6\cdot 4^2)$ ($Z_{11}$), $\Delta(6\cdot 6^2)$, $\Delta(6\cdot 6^2)$ ($Z_{1,0}$), and $\Delta(6\cdot 3^2)$ which corresponds to the elliptic complete intersection singularity $\delta 1$ in C.~T.~C.~Wall's notation \cite{Wa2}. (For a list of Fuchsian hypersurface and complete intersection singularities see \cite{Ebeling03}.) These are 6 cases.
\begin{table}
\begin{center}
\begin{tabular}{|c|c|c|c|c|c|}
\hline
 $G$    & $|G|$ & $w,x,y,z(,u)$ & $c_G$ & $R(0,x,y,z(,u))$ & Sing.\\
\hline
$\Delta(3 \cdot 2^2)$  & $12$ & $2,3,4,6$ & 1& $z^2+4y^3 +27x^4$ & $E_6$ \\
$\Delta(3 \cdot 3^2)$ & $27$ & $3,3,6,9$ & 3 & $z^2+4y^3+27x^6$ & $\widetilde{E}_8$ \\
$\Delta(3 \cdot 4^2)$ & $48$ & $4,3,8,12$ & 1 & $z^2+4y^3+27x^8$ & $E_{14}$ \\
$\Delta(3 \cdot 5^2)$ & $75$ & $5,3,10,15;30$ & 1 & $z^2+4y^3+27x^{10}$ & $E_{18}$ \\
$\Delta(6\cdot 2^2)$ & $24$ & $2,4,6,9$ & 1 & $z^2+4xy^3+27x^3$ & $E_7$ \\
$\Delta(6\cdot 3^2)$ & $54$ & $6,6,6,6,9$ & 3 & $\left\{ \begin{array}{c} z^2-xy \\ u^2 +4xyz+27x^3 \end{array} \right\}$ & $\delta 1$ \\
$\Delta(6\cdot 4^2)$ & $96$ & $4,6,8,15$ & 1 &  $z^2+4xy^3+27x^5$ & $Z_{11}$ \\
$\Delta(6 \cdot 5^2)$ & $150$ & $10,6,8,10,15$ & 1 & $\left\{ \begin{array}{c} z^2-xy \\ u^2 +4x^2yz+27x^5 \end{array} \right\}$ & no name \\
$\Delta(6\cdot 6^2)$ & $216$ & $6,6,12,21$ & 3 & $z^2+4xy^3+27x^7$ & $Z_{1,0}$ \\
$\Delta(6 \cdot 7^2)$ & $294$ & $14,6,10,14,21$ & 1 & $\left\{ \begin{array}{c} z^2-xy \\ u^2 +4x^3yz+27x^7 \end{array} \right\}$ & no name \\
$\Delta(6\cdot 8^2)$ & $384$ & $8,6,16,27;54$ & 1 & $z^2+4xy^3+27x^9$ & $Z_{19}$ \\
\hline
\end{tabular}
\end{center}
\caption{The first subgroups of types (C) and (D) and surface singularities} \label{TabCD}
\end{table}
The remaining 8 exceptional subgroups of types (E)--(L) all correspond to Fuchsian singularities except in the case (H) which is the icosahedral group $I$ corresponding to the Kleinian singularity $E_8$. 
Altogether we have 14 cases which we will consider in this paper. They are listed in Table~\ref{TabSL3sing}.
\begin{table}
\begin{center}
\begin{tabular}{|c|c|c|c|c|c|}
\hline
 $G$    & $|G|$ & $w,x,y,z(,u)$ & $c_G$ & $R(0,x,y,z(,u))$ & Sing.\\
\hline
(C): $T$ & $12$ & $2,3,4,6$ & 1& $z^2+4y^3 +27x^4$ & $E_6$ \\
$\Delta(3 \cdot 4^2)$ & $48$ & $4,3,8,12$ & 1 & $z^2+4y^3+27x^8$ & $E_{14}$ \\
(D): $O$ & $24$ & $2,4,6,9$ & 1 & $z^2+4xy^3+27x^3$ & $E_7$ \\
$\Delta(6\cdot 3^2)$ & $54$ & $6,6,6,6,9$ & 3 & $\left\{ \begin{array}{c} z^2-xy \\ u^2 +4xyz+27x^3 \end{array} \right\}$ & $\delta 1$ \\
$\Delta(6\cdot 4^2)$ & $96$ & $4,6,8,15$ & 1 &  $z^2+4xy^3+27x^5$ & $Z_{11}$ \\
$\Delta(6\cdot 6^2)$ & $216$ & $6,6,12,21$ & 3 & $z^2+4xy^3+27x^7$ & $Z_{1,0}$ \\
(E) &108 & $6,6,9,12,12$ & 3 & $\left\{ \begin{array}{c} 9u^2-12z^2 \\432y^2-x^3-36xz \end{array} \right\}$ & $K'_{1,0}$  \\
 (F) & 216 & $6,9,12,12$ & 3 & $\begin{array}{l} 4z^3-144yz^2\\ +1728y^2z -186624x^4 \end{array}$ & $U_{12}$\\
 (G) & 648 & $9,12,18,18$ & 6 &  $4z^3-9yz^2+6y^2z-y^3+6912x^3y$ & $U_{1,0}$ \\
(H)=$I$ & $60$ & $2,6,10,15$ & 1 & $z^2-y^3+1728x^5$ & $E_8$\\
 (I)  & $168$ & $4,6,14,21$ & 1 & $z^2-y^3-1728x^7$ & $E_{12}$ \\
 (J) & 180 & $6,6,12,15$ & 3 & $y^3-xz^2+64x^2y^2$ & $Q_{2,0}$  \\
 (K) & 504 & $6,12,18,21$ & 3 & $y^3-xz^2-256x^3y$ & $Q_{11}$ \\
 (L) & 1080 & $6,12,30,45$ & 3 & $\begin{array}{l} 459165024z^2-25509168y^3 \\-(7558272-2519424\sqrt{15}i)x^5y \end{array}$ & $E_{13}$\\
\hline
\end{tabular}
\end{center}
\caption{Subgroups of ${\rm SL}_3(\CC)$ and surface singularities} \label{TabSL3sing}
\end{table}
These singularities are surface singularities and they are isolated except in the three cases $\Delta(6\cdot 3^2)$, (E) and (J). They correspond to Kleinian singularities in the cases $T$, $O$ and (H) (the icosahedral group $I$) and to Fuchsian singularities in the other cases. They correspond to simple ($T$, $O$, $I$), unimodal ($\Delta(3 \cdot 4^2)$, $\Delta(6\cdot 4^2)$, (F), (I), (K), (L)) and bimodal ($\Delta(6\cdot 6^2)$, (G), (J)) hypersurface singularities, to the unimodal complete intersection singularity of type $K'_{1,0}$ (type (E)) in Wall's notation  \cite{Wa1}, and to the elliptic complete intersection singularity $\delta1$. The names of the hypersurface singularities according to V.~I.~Arnold's classification \cite{Ar} are indicated in the last column of Table~\ref{TabSL3sing}.

\section{Poincar\'e series and Coxeter elements}

\begin{sloppypar}
We now consider the isolated singularities corresponding to these singularities. They are quasihomogeneous. This means the following.  A complex polynomial 
$f(x_1, \ldots, x_n)$ is called quasihomogeneous, if there are positive integers $w_1, \ldots, w_n$ (called {\em weights}) and $d$ (called {\em degree}) such that $f(\lambda^{w_1}x_1, \ldots , \lambda^{w_n}x_n) = \lambda^d f(x_1, \ldots , x_n)$ for $\lambda \in \CC^\ast$. A complete intersection singularity given as the zero set of polynomials $f_1(x_1, \ldots,x_n), \ldots , f_k(x_1, \ldots , x_n)$ is called quasihomogeneous, if $f_1, \ldots, f_k$ are quasihomogeneous with respect to the same weights $w_1, \ldots, w_n$ but degrees $d_1, \ldots, d_k$ respectively. We call the system $W:=(w_1, \ldots, w_n;d_1, \ldots , d_k)$ the {\em weight system} corresponding to the set of polynomials. Let $c_W$ be the greatest common divisor of $w_1, \ldots, w_n,d_1, \ldots , d_k$. The weight system is called {\em reduced} if $c_W=1$. 

We assume that $f_1(0)= \cdots = f_k(0)=0$ and the system of equations $f_1= \cdots =f_k =0$ has an isolated singularity at the origin. The  coordinate algebra $A_f := \CC[x_1, \ldots , x_n]/(f_1, \ldots, f_k)$ is a $\ZZ$-graded algebra with respect to the system of weights $(w_1, \ldots, w_n; d_1, \ldots, d_k)$. Therefore we can consider the decomposition of $A_f$ as a $\ZZ$-graded $\CC$-vector space:
\[ A_f:= \bigoplus_{k=0}^\infty A_{f,k}, \quad A_{f,k}:= \left\{ g \in A_f \, \left| \, g(\lambda^{w_1}x_1, \ldots , \lambda^{w_n}x_n) = \lambda^k g(x_1, \ldots , x_n) \right\} \right. .
\]
The formal power series $p_f(t):= \sum_{k=0}^\infty (\dim_\CC A_{f.k})t^k$ is called the {\em Poincar\'e series} of $A_f$. It is given by
\[ p_f(t)= \frac{\prod_{i=1}^k (1-t^{d_i})}{\prod_{j=1}^n (1-t^{w_j})} \, .
\]

\end{sloppypar}

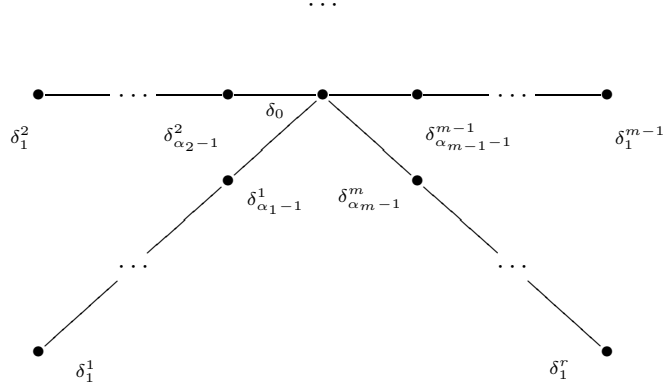
\begin{figure}
$$
\xymatrix{ & & & {\cdots} & & & \\
 *{\bullet} \ar@{-}[r] \ar@{}_{\delta^2_1}[d]  & {\cdots} \ar@{-}[r]  & *{\bullet} \ar@{-}[r]   \ar@{}_{\delta^2_{\alpha_2-1}}[d] & *{\bullet} \ar@{-}[dl]  \ar@{-}[dr] \ar@{-}[r]  \ar@{}^{\delta_{0}}[l]  & *{\bullet} \ar@{-}[r]  \ar@{}^{\delta^{m-1}_{\alpha_{m-1}-1}}[d]  & {\cdots} \ar@{-}[r]  &*{\bullet} \ar@{}^{\delta^{m-1}_1}[d]   \\
& &  *{\bullet} \ar@{-}[dl] \ar@{}_{\delta^1_{\alpha_1-1}}[r]  &  & *{\bullet} \ar@{-}[dr] \ar@{}^{\delta^m_{\alpha_m-1}}[l] &  &  \\
 & {\cdots} \ar@{-}[dl] & &  & & {\cdots} \ar@{-}[dr] & \\
*{\bullet}  \ar@{}_{\delta^1_1}[r] & & &  & & & *{\bullet}  
\ar@{}^{\delta^r_1}[l]
  } 
$$
\caption{The graph $T^-_{\alpha_1,\ldots ,\alpha_m}$} \label{FigT-pqr}
\end{figure}

Let $(X,x)$ be a Kleinian singularity. Then the minimal resolution of the singularity $x$ has an exceptional divisor with the dual graph depicted in Fig.~\ref{FigT-pqr} with $m=1$ in the case of the $A_n$-singularities and $m=3$ in the other cases (see, e.g., \cite{Br}). Here all vertices correspond to rational curves of self-intersection number $-2$, the mutual intersection numbers are either 0 or 1, and two vertices are joined by an edge if and only if the intersection number of the corresponding rational curves is equal to 1. The values of the numbers $\alpha_1, \ldots , \alpha_m$ are indicated in Table~\ref{TabSL2sing}. They are the {\em Dolgachev numbers} of the singularity, see \cite{ET}. It turns out that these graphs are precisely the ordinary Coxeter-Dynkin diagrams of type ADE. (Note that the corresponding intersection matrix is the Cartan matrix multiplied by $-1$.)

Now let $(X,x)$ be a Fuchsian hypersurface or complete intersection singularity. A natural compactification of $X$ is given by $\overline{X}:= {\rm Proj}(A_f[t])$, where $t$ has degree 1 for the grading of $A_f[t]$ (see \cite{Pi}). This is a normal projective surface with a $\CC^\ast$-action. The surface $\overline{X}$ may acquire additional singularities on the boundary $\overline{X}_\infty:= \overline{X} \setminus X = {\rm Proj}(A_f)$. Let $g=g(\overline{X}_\infty)$ be the genus of the boundary. We assume $g=0$. Let $\pi: S \to \overline{X}$ be the minimal normal crossing resolution of all singularities of $\overline{X}$. The preimage $\widetilde{X}_\infty$ of $\overline{X}_\infty$ under $\pi \colon S \to \overline{X}$
consists of the strict transform $\delta_0$ of $\overline{X}_\infty$ and $m$ chains
$\delta^i_1, \ldots , \delta^i_{\alpha_i-1}$, $i=1, \ldots , m$, of rational curves of
self-intersection $-2$ which intersect again according to the dual graph shown in
Figure~\ref{FigT-pqr} (see, e.g., \cite{Dolgachev07,Ebeling03}). By the adjunction formula and $g=0$, the self-intersection number of 
the rational curve $\delta_0$ is also $-2$. The numbers $\alpha_1, \ldots , \alpha_m$ of the Fuchsian singularities corresponding to finite subgroups of ${\rm SL}_3(\CC)$ are indicated in Table~\ref{TabPoincE}. They are again the {\em Dolgachev numbers} of the singularity, see \cite{ET, Ebeling99}.
\begin{table}
\begin{center}
\begin{tabular}{|c|c|c|c|c|}
\hline
$G$  & Name & Normal form    & Weights & $\alpha_1, \ldots , \alpha_m$\\
\hline
(C): $T$  &  $E_6$ & $z^2 + y^3 + x^4$ & 3,4,6;12 & ${2,3,3}$ \\
$\Delta(3 \cdot 4^2)$ & $E_{14}$ & $z^2+y^3+x^8$ & 3,8,12;24 & ${3,3,4}$ \\
(D): $O$  & $E_7$ & $z^2+y^3+yx^3$ & 4,6,9;18 & ${2,3,4}$ \\
$\Delta(6\cdot 3^2)$ & $\delta 1$ & $\left\{ \begin{array}{c} xy+z^2 \\ x^3+y^3+z^3+w^2 \end{array} \right\}$ & 2,2,2,3;4,6 & ${2,2,2,2,2,2}$ \\
$\Delta(6\cdot 4^2)$ & $Z_{11}$ & $z^2+xy^3+x^5$ & 6,8,15;30 & ${2,3,8}$ \\
$\Delta(6\cdot 6^2)$ & $Z_{1,0}$ & $z^2+xy^3+x^7$ & 2,4,7;14 & ${2,2,2,4}$ \\
(E)  & $K'_{1,0}$ & $\left\{ \begin{array}{c} xu+y^2 \\ax^4+xy^2+z^2+u^2, \\ a \neq 0, \frac{1}{4} \end{array} \right\}$   & 2,3,4,4;6,8 & ${2,2,4,4}$ \\
(F)  & $U_{12}$ & $z^3 +y^3+x^4$ & 3,4,4;12& ${4,4,4}$ \\
(G)  & $U_{1,0}$ & $z^3+yz^2+x^3y$ & 2,3,3;9 & ${2,3,3,3}$ \\
(H)=$I$  & $E_8$ & $z^2+y^3+x^5$ & 6,10,15;30 & ${2,3,5}$ \\
(I)  & $E_{12}$  & $z^2+y^3+x^7$ & 6,14,21;42 & ${2,3,7}$  \\
(J)  &   $Q_{2,0}$ & $xz^2+y^3+x^4y$ & 2,4,5;12 & ${2,2,2,5}$ \\
(K) &  $Q_{11}$ & $xz^2+y^3+yx^3$ & 4,6,7;18 & ${2,4,7}$ \\
(L)  & $E_{13}$ & $z^2+y^3+x^5y$ & 4,10,15;30 & ${2,4,5}$ \\
\hline
\end{tabular}
\end{center}
\caption{Normal forms, reduced weight systems, and Dolgachev numbers} \label{TabPoincE}
\end{table}

We call the graph $T^-_{\alpha_1,\ldots ,\alpha_m}$ a {\em Coxeter-Dynkin diagram}. Let $V_-$ be the free $\ZZ$-module with the basis
\[ \delta^1_1, \ldots, \delta^1_{\alpha_1-1}, \delta^2_1, \ldots ,\delta^2_{\alpha_2-1}, \ldots , \delta^m_1, \ldots, \delta^m_{\alpha_m-1}, \delta_0
\]
equipped with the symmetric bilinear form $\langle-,-\rangle$ given by the intersection matrix corresponding to Fig.~\ref{FigT-pqr}. This defines a lattice $(V_-,\langle-,-\rangle)$.
\begin{figure}
$$
\xymatrix{ & & & {\cdots} & & & \\
 & & & *{\bullet}  \ar@{==}[d] \ar@{-}[dr]  \ar@{-}[ldd] \ar@{-}[rdd] \ar@{}_{\delta_1}[l]
 & & &  \\
 *{\bullet} \ar@{-}[r] \ar@{}_{\delta^2_1}[d]  & {\cdots} \ar@{-}[r]  & *{\bullet} \ar@{-}[r] \ar@{-}[ur]   \ar@{}_{\delta^2_{\alpha_2-1}}[d] & *{\bullet} \ar@{-}[dl]  \ar@{-}[dr] \ar@{-}[r]  \ar@{}^{\delta_{0}}[l]  & *{\bullet} \ar@{-}[r]  \ar@{}^{\delta^{m-1}_{\alpha_{m-1}-1}}[d]  & {\cdots} \ar@{-}[r]  &*{\bullet} \ar@{}^{\delta^{m-1}_1}[d]   \\
& &  *{\bullet} \ar@{-}[dl] \ar@{}_{\delta^1_{\alpha_1-1}}[r]  &  & *{\bullet} \ar@{-}[dr] \ar@{}^{\delta^m_{\alpha_m-1}}[l] &  &  \\
 & {\cdots} \ar@{-}[dl] & &  & & {\cdots} \ar@{-}[dr] & \\
*{\bullet}  \ar@{}_{\delta^1_1}[r] & & &  & & & *{\bullet}  
\ar@{}^{\delta^r_1}[l]
  } 
$$
\caption{The graph $T_{\alpha_1,\ldots ,\alpha_m}$} \label{FigTpqr}
\end{figure}
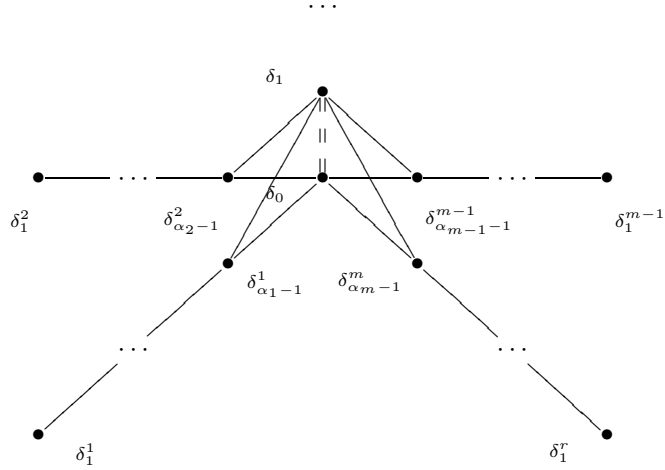

\begin{figure}
$$
\xymatrix{ & & & {\cdots} & & & \\
& & & *{\bullet} \ar@{-}[d] \ar@{}_{\delta_2}[l] & & & \\
 & & & *{\bullet}  \ar@{==}[d] \ar@{-}[dr]  \ar@{-}[ldd] \ar@{-}[rdd] \ar@{}_{\delta_1}[l]
 & & &  \\
 *{\bullet} \ar@{-}[r] \ar@{}_{\delta^2_1}[d]  & {\cdots} \ar@{-}[r]  & *{\bullet} \ar@{-}[r] \ar@{-}[ur]   \ar@{}_{\delta^2_{\alpha_2-1}}[d] & *{\bullet} \ar@{-}[dl]  \ar@{-}[dr] \ar@{-}[r]  \ar@{}^{\delta_{0}}[l]  & *{\bullet} \ar@{-}[r]  \ar@{}^{\delta^{m-1}_{\alpha_{m-1}-1}}[d]  & {\cdots} \ar@{-}[r]  &*{\bullet} \ar@{}^{\delta^{m-1}_1}[d]   \\
& &  *{\bullet} \ar@{-}[dl] \ar@{}_{\delta^1_{\alpha_1-1}}[r]  &  & *{\bullet} \ar@{-}[dr] \ar@{}^{\delta^m_{\alpha_m-1}}[l] &  &  \\
 & {\cdots} \ar@{-}[dl] & &  & & {\cdots} \ar@{-}[dr] & \\
*{\bullet}  \ar@{}_{\delta^1_1}[r] & & &  & & & *{\bullet}  
\ar@{}^{\delta^r_1}[l]
  } 
$$
\caption{The graph $T^+_{\alpha_1,\ldots ,\alpha_m}$} \label{FigT+pqr}
\end{figure}
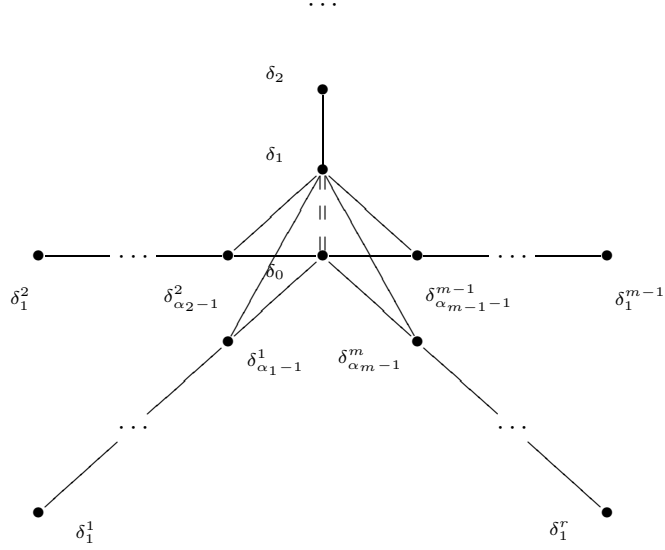

We consider two extensions of this lattice. Let $V_0= V_- \oplus \ZZ \delta_1$ with the symmetric bilinear form defined by Fig.~\ref{FigTpqr}. Here the double dashed line between $\delta_0$ and $\delta_1$ means $\langle \delta_0, \delta_1 \rangle =-2$. Let $V_+=V_0 \oplus \ZZ \delta_2$ with the symmetric bilinear form defined by Fig.~\ref{FigT+pqr}.

If $(V,\langle-,-\rangle)$ is an arbitrary lattice and $e \in V$ is a \emph{root},
i.e.\ $\langle e,e \rangle=-2$, then the reflection corresponding to $e$ is 
defined by 
\[ s_e(x) = x - \frac{2\langle x,e \rangle}{\langle e, e \rangle} e 
          = x + \langle x,e\rangle e \quad \text{for } x \in V. \]
If $B=(e_1, \ldots , e_n)$ is an ordered basis consisting of roots, then the 
\emph{Coxeter element} $\tau$ corresponding to $B$ is defined by
\[ \tau= s_{e_1} s_{e_2} \cdots s_{e_n}. \] 
For a Coxeter element $\tau$, let 
$\Delta(t)= \det (1 - \tau^{-1}t )$ be its characteristic polynomial, using a suitable normalization.

If $D$ is a Coxeter-Dynkin diagram, then we denote by $\Delta_D(t)$ the characteristic polynomial of the Coxeter element corresponding to the graph $D$. These polynomials can be computed as in \cite{Ebeling83} and one gets
\begin{eqnarray*}
\Delta_{T^-_{\alpha_1, \ldots , \alpha_m}}(t) & = & (1+t) \prod_{i=1}^m \frac{1-t^{\alpha_i}}{1-t} - t \sum_{i=1}^m \frac{1-t^{\alpha_i-1}}{1-t} \prod_{j=1 \atop j \neq i}^m \frac{1-t^{\alpha_j}}{1-t}, \\
\Delta_{T_{\alpha_1, \ldots , \alpha_m}}(t) & = & (1-t)^{2-m}(1-t^{\alpha_1}) \cdots (1-t^{\alpha_m}), \\
\Delta_{T^+_{\alpha_1, \ldots , \alpha_m}}(t) & = & (1-2t-2t^2+t^3) \prod_{i=1}^m \frac{1-t^{\alpha_i}}{1-t}+t^2 \sum_{i=1}^m \frac{1-t^{\alpha_i-1}}{1-t} \prod_{j=1 \atop j \neq i}^m \frac{1-t^{\alpha_j}}{1-t}.
\end{eqnarray*}
(The last two formulas can also be found in \cite[p.~98]{Ebeling87}, but note that, unfortunately, there is a misprint in \cite[p.~98]{Ebeling87}.)

Now we can state the main result of \cite{EP1}.

\begin{theorem} \label{theo:sing}
\begin{itemize}
\item[{\rm (i)}] For a Kleinian singularity not of type $A_{2n}$ we have
\[ p_f(t) = \frac{\Delta_{T^-_{\alpha_1, \ldots , \alpha_m}}(t)}{\Delta_{T_{\alpha_1, \ldots , \alpha_m}}(t)}.
\]
\item[{\rm (ii)}] For a Fuchsian singularity with $g=0$ we have
\[ p_f(t) = \frac{\Delta_{T^+_{\alpha_1, \ldots , \alpha_m}}(t)}{\Delta_{T_{\alpha_1, \ldots , \alpha_m}}(t)}.
\]
\end{itemize}
\end{theorem}

\begin{remark} Unfortunately, the exclusion of the case $A_{2n}$ is only implicit in \cite{EP1} and was forgotten in the statement of \cite[Theorem~1]{EP1}. 
\end{remark} 

\begin{remark}
Note that we have $T_{2,3,3} \sim T^-_{3,3,3}$, $T_{2,3,4} \sim T^-_{2,4,4}$, $T_{2,3,5} \sim T^-_{2,3,6}$, where $\sim$ means equivalence under the braid group action, see \cite{Ebeling16}. Similarly, one can show that the graphs $T_{2n-1}$, $n \geq 1$, and $T_{2,2,n}$, $n \geq 2$, are equivalent under the braid group action to the extended Coxeter-Dynkin diagrams of type $A_{2n-1}$ and $D_{n+2}$ respectively.
\end{remark}

\section{Poincar\'e series of subgroups of ${\rm SL}_2(\CC)$ and ${\rm SL}_3(\CC)$}
Let $G$ be a finite subgroup of ${\rm SL}_n(\CC)$ for $n=2,3$.
Consider the algebra of complex polynomials $\CC[x_1, \ldots , x_n]$ graded by the degree for homogeneous ones. It is isomorphic to the symmetric algebra
\[ S:=S(\CC^n) = \bigoplus_{k=0}^\infty S^k(\CC^n), \]
where $S^k(\CC^n)$ denotes the $k$-th symmetric power of $\CC^n$. Let $S^G$ be the algebra of invariant polynomials with respect to $G$. 

For $n=2$, it is generated by 3 elements $x,y,z$ which satisfy a relation $R(x,y,z)=0$. The elements $x,y,z$ correspond to invariant polynomials and their degrees correspond to the weights of these variables. Let $c_G$  denote the greatest common divisor of these weights. 
The weights of the variables $x,y,z$, the number $c_G$, and the polynomial $R(x,y,z)$ are indicated in Table~\ref{TabSL2sing}. 

Now let $G$ be one of the finite subgroups of ${\rm SL}_3(\CC)$ of Table~\ref{TabSL3sing}. Except in the cases (E) and $\Delta(6\cdot 3^2)$, the algebra $S^G$ is generated by 4 elements $w,x,y,z$ which satisfy a relation $R(w,x,y,z)$.  In the cases (E) and $\Delta(6\cdot 3^2)$, $S^G$ is generated by 5 elements $w,x,y,z,u$ which satisfy two relations $R_1(w,x,y,z,u)=0$ and $R_2(w,x,y,z,u)=0$. The degrees of the invariants and the polynomials $R(w,x,y,z)$ and $R_1(w,x,y,z,u)$, $R_2(w,x,y,z,u)$ respectively can be found in \cite{YY}.  The degrees of the invariant polynomials and the restriction to the hyperplane $w=0$ of the polynomials $R(w,x,y,z)$ and $R_1(w,x,y,z,u)$, $R_2(w,x,y,z,u)$ respectively are indicated in Table~\ref{TabSL3sing}. Let $c_G$ be the greatest common divisor of the weights of the remaining variables $x,y,z(,u)$ (with the weight of $w$ excluded). The number $c_G$ is also indicated in Table~\ref{TabSL3sing}. Note that, except in the case (G), $c_G$ also divides the weight of $w$.

For $n=2$, the algebra $A_G:=S^G=\CC[x,y,z]/R(x,y,z)$ coincides with the coordinate algebra $A_f$ of the corresponding singularity indicated in the last column of Table~\ref{TabSL2sing} up to the grading. The grading is shifted by $c_G$. For $n=3$ and $G$ one of the cases of Table~\ref{TabSL3sing} except the cases (E) and $\Delta(6\cdot 3^2)$, the algebra $A_G:=\CC[x,y,z]/R(0,x,y,z)$ coincides with the coordinate algebra $A_f$ of the corresponding singularity indicated in the last column of Table~\ref{TabSL3sing} with the grading shifted by $c_G$. In the cases (E) and $\Delta(6\cdot 3^2)$, the algebra $A_G:=\CC[x,y,z,u]/(R_1(0,x,y,z,u),R_2(0,x,y,z,u))$ coincides with the coordinate algebra $A_f$ of the complete intersection singularity $K'_{1,0}$ and $\delta 1$ respectively, again with the grading shifted by $c_G$. Let $p_G(t)$ be the Poincar\'e series of the algebra of $A_G$. Then we have
\[ p_G(t) = p_f(t^{c_G}) \mbox{ for } G \subset {\rm SL}_2(\CC), \quad p_G(t)=\frac{p_f(t^{c_G})}{(1-t^{\deg w})} \mbox{ for } G \subset {\rm SL}_3(\CC).
\]

The finite subgroups $G \subset {\rm SL}_n(\CC)$ for $n=2,3$ under consideration have a natural $n$-dimensional representation $\gamma$ which is irreducible (except in the cases $G=\mathcal{C}_l$). Let $\gamma^\ast$ be its contragredient representation.  Let $\gamma_0, \ldots , \gamma_l$ be the irreducible representations of $G$, where $\gamma_0$ is the trivial representation.  Let $B=(b_{ij})$  and $B^\ast=(b^\ast_{ij})$ be the $(l+1) \times (l+1)$-matrices defined by decomposing the tensor products 
\[\gamma_j \otimes \gamma = \bigoplus_i b_{ij} \gamma_i \quad \mbox{and} \quad \gamma_j \otimes \gamma^\ast = \bigoplus_i b^\ast_{ij} \gamma_i
\] 
respectively into irreducible components. 

For each integer $k \geq 0$, let $\rho_k$ be the representation of $G$ on $S^k(\CC^n)$ induced by its natural action on $\CC^n$. We have a decomposition $\rho_k= \sum_{i=0}^l v_{ki} \gamma_i$ with $v_{ki} \in \ZZ$. We associate to $\rho_k$ the vector $v_m=(v_{m0}, \ldots , v_{ml})^t \in \ZZ^{l+1}$.
As in \cite[p.~211]{Kostant85} we define
$$P_G(t):= \sum_{m=0}^\infty v_m t^m.$$
This is a formal power series with coefficients in $\ZZ^{l+1}$. We also put 
$P_G(t)_i:= \sum_{m=0}^\infty v_{mi} t^m$. Note that $P_G(t)_0$ is the usual Poincar\'e series $p_G(t)$ of the group $G$. Let $V$ denote the set of all formal power series 
$x = \sum_{m=0}^\infty x_m t^m$
with $x_m \in \ZZ^{l+1}$. This is a free module of rank $l+1$ over the ring $R$ of
formal power series with integer coefficients. 

Now let $n=2$ and $G \subset {\rm SL}_2(\CC)$ be a finite subgroup not of type $\mathcal{C}_{2n+1}$. Then $c_G=2$ and we have
\[
p_f(t^2)=P_G(t)_0.
\]
Moreover, we have $\gamma^\ast=\gamma$ and therefore $B^\ast=B$. The irreducible representations of ${\rm SL}_2(\CC)$ are of the form $\rho_m$, $m$ a non-negative integer. The Clebsch-Gordon formula reads in this case
\[  \rho_m \otimes \gamma = \rho_{m+1} \oplus \rho_{m-1} \]
setting $\rho_{-1}=0$ (cf., e.g., \cite[Exercise~11.11]{FH}). This yields the equation
\[ Bv_m=v_{m+1}+v_{m-1}.
\]
Following \cite[p.~222]{Kostant85},  one can easily derive from this equation that $x=P_G(t)$ is a solution of the following linear equation in $V$:
\[ ((1+t^2)I -tB)x=v_0.
\]
Let $M(t)$ be the matrix $(1+t^2)I -tB$ and $M_0(t)$ be the matrix obtained by replacing the first column of $M(t)$ by $v_0=(1,0, \ldots, 0)^t$. By Cramer's rule we can derive the following theorem \cite[Sect.~3]{Ebeling02} (see also \cite{St}).
\begin{theorem}
For a finite subgroup $G \subset {\rm SL}_2(\CC)$ not of type $\mathcal{C}_{2n+1}$ we have
\[
p_f(t^2)=P_G(t)_0 = \frac{\det M_0(t)}{\det M(t)} = \frac{\det(t^2I -\tau)}{\det(t^2I-\tau_a)},
\]
where $\tau$ is the Coxeter element and $\tau_a$ the affine Coxeter element of the corresponding root system of type ADE associated to the singularity defined by the equation $f=0$ with the same name.
\end{theorem}

Now let $n=3$ and $G \subset  {\rm SL}_3(\CC)$ be a finite subgroup. For a pair $a,b$ of non-negative integers, let $\Gamma_{a,b}$ be the unique irreducible, finite-dimensional representation of ${\rm SL}_3(\CC)$ of \cite[Theorem~13.1]{FH}. By \cite[Proposition~15.25]{FH} and \cite[(13.5)]{FH}, we have for a non-negative integer $m$ (setting $\Gamma_{-1,b} =0$) the following Clebsch-Gordon formulas:
\begin{eqnarray*}
\Gamma_{m,0} \otimes \gamma & = & \Gamma_{m+1,0} \oplus \Gamma_{m-1,1}, \\
\Gamma_{m,0} \otimes \gamma^\ast & = & \Gamma_{m-1,0} \oplus \Gamma_{m,1}.
\end{eqnarray*}
Since $\Gamma_{m,0} = \rho_m$, we can derive from these formulas
\[ 
v_{m+2}=Bv_{m+1} - B^\ast v_m+ v_{m-1}.
\]
Therefore $x=P_G(t)$ is a solution of the following linear
equation in $V$ (see also \cite{BI,BP}):
\[
((1-t^3)I-tB+t^2B^\ast)x=v_0.
\]
Let $M(t)$ be the matrix $(1-t^3)I-tB+t^2B^\ast$ and $M_0(t)$ be the matrix obtained by replacing the first column of $M(t)$ by $v_0=(1,0, \ldots, 0)^t$. Again Cramer's rule yields
\begin{equation*} \label{eq:P}
P_G(t)_0= \frac{\det M_0(t)}{\det M(t)}.
\end{equation*}

We have the following theorem:

\begin{theorem} \label{thm:main}
Let $G \subset  {\rm SL}_3(\CC)$ be one of the groups $T$, $\Delta(3 \cdot 4^2)$, $O$, $\Delta(6 \cdot 3^2)$, $\Delta(6 \cdot 4^2)$, $\Delta(6 \cdot 6^2)$, (E), (F), (G), (H)=$I$, (I), (J), (K), or (L), let $c_G$ be the greatest common divisor of the weights of the variables $x,y,z(,u)$, and let $\alpha_1, \ldots, \alpha_m$ be the Dolgachev numbers of the singularity corresponding to $G$. Moreover, let $q^{(e)}_{a,b}(t)=(1-t)^a(1-t^e)^b$ for $a,b,e \in \ZZ$. 
\begin{itemize}
\item[{\rm (i)}] For $G=T,O,I$ ($E_6, E_7, E_8$) we have $c_G=1$ and
\[ \det M_0(t) = q^{(2)}_{a,b}(t) \Delta_{T^-_{\alpha_1, \alpha_2, \alpha_3}}(t), \quad \det M(t)= (1-t) q^{(2)}_{a,b}(t) \Delta_{T_{2,\alpha_1, \alpha_2, \alpha_3}}(t),
\]
where $(a,b)=(3,0), (3,1), (4,0)$ respectively.
\item[{\rm (ii)}]
For $G$=(I), $\Delta(3 \cdot 4^2)$, $\Delta(6 \cdot 4^2)$ ($E_{12}, E_{14}, Z_{11}$) we have $c_G=1$ and 
\[ \det M_0(t) = q^{(4)}_{a,b} (t)\Delta_{T^+_{\alpha_1, \alpha_2, \alpha_3}}(t), \quad \det M(t)= (1-t) q^{(4)}_{a,b}(t) \Delta_{T_{4, \alpha_1, \alpha_2, \alpha_3}}(t).
\]
where $(a,b)=(3,0), (3,2), (2,3)$ respectively.
\item[{\rm (iii)}] 
For $G$=(K), (L), (F), $\Delta(6 \cdot 6^2)$, (J), (E)  ($Q_{11}, E_{13}, U_{12}, Z_{1,0}, Q_{2,0}, K'_{1,0}$) we have $c_G=3$ and
\[ \det M_0(t) = q^{(2)}_{a,b}(t^3) \Delta_{T^+_{\alpha_1, \ldots, \alpha_m}}(t^3), \quad \det M(t)= (1-t^3) q^{(2)}_{a,b}(t^3) \Delta_{T_{2,\alpha_1, \ldots, \alpha_m}}(t^3),
\]
where $(a,b)=(6,-1), (7,-1), (7,-2), (7,1), (8,-2), (8,-3)$ respectively.
\item[{\rm (iv)}] For $G=\Delta(6 \cdot 3^2)$  ($\delta1$) we have $c_G=3$, $m=6$, $\alpha_i=2$ for $i=1, \ldots ,m$, and
\[ \det M_0(t) = (1-t^3)q^{(2)}_{9,-3}(t^3) \Delta_{T^+_{2,\alpha_2, \ldots , \alpha_m}}(t^3), \quad \det M(t)= q^{(2)}_{9,-3}(t^3) \Delta_{T_{\alpha_2, \ldots , \alpha_m}}(t^3).
\]
\item[{\rm (v)}] 
For $G$=(G) ($U_{1,0}$) we have $c_G=6$ and 
\[ \det M_0(t) = q(t^3) \Delta_{T^+_{\alpha_1,\alpha_2,\alpha_3, \alpha_4}}(t^6), \quad \det M(t)= (1-t^9)q(t^3) \Delta_{T_{\alpha_1,\alpha_2,\alpha_3, \alpha_4}}(t^6),
\]
where $q(t)=\frac{(1-t)^4(1-t^2)(1-t^3)}{(1-t^6)}$. 
\end{itemize}
\end{theorem}

\begin{table}
\begin{center}
\begin{tabular}{|c|c|c|}
\hline
 $G$    & $\det M_0(t)$ & $\det M(t)$ \\
\hline
 $T$  & $(1-t)^3\Delta_{T^-_{2,3,3}}(t)$ & $(1-t)^4 \Delta_{T_{2,2,3,3}}(t)$ \\
$\Delta(3 \cdot 4^2)$ & $(1-t)^3(1-t^4)^2 \Delta_{T^+_{3,3,4}}$ & $(1-t)^4  (1-t^4)^2\Delta_{T_{3,3,4,4}}(t)$\\
$O$  & $(1-t)^3 (1-t^2) \Delta_{T^-_{2,3,4}}(t)$ & $(1-t)^4 (1-t^2) \Delta_{T_{2,2,3,4}}(t)$  \\
 $\Delta(6\cdot 3^2)$ & $\frac{(1-t^3)^{10}}{(1-t^6)^3} \Delta_{T^+_{2,2,2,2,2,2}}(t^3)$ & $\frac{(1-t^3)^9}{(1-t^6)^3} \Delta_{T_{2,2,2,2,2}}(t^3)$ \\
 $\Delta(6\cdot 4^2)$ & $(1-t)^2(1-t^4)^3 \Delta_{T^+_{2,3,8}}(t)$ & $(1-t)^3(1-t^4)^3 \Delta_{T_{2,3,4,8}}(t)$ \\
 $\Delta(6\cdot 6^2)$ & $(1-t^3)^7(1-t^6) \Delta_{T^+_{2,2,2,4}}(t^3)$ & $(1-t^3)^8(1-t^6) \Delta_{T_{2,2,2,2,4}}(t^3)$ \\
 (E)  & $\frac{(1-t^3)^8}{(1-t^6)^3} \Delta_{T^+_{2,2,4,4}}(t^3)$ & $\frac{(1-t^3)^9}{(1-t^6)^3} \Delta_{T_{2,2,2,4,4}}(t^3)$  \\
 (F)  & $\frac{(1-t^3)^7}{(1-t^6)^2} \Delta_{T^+_{4,4,4}}(t^3)$ & $\frac{(1-t^3)^8}{(1-t^6)^2} \Delta_{T_{2,4,4,4}}(t^3)$  \\
   (G)  & $\frac{(1-t^3)^4(1-t^6)(1-t^9)}{(1-t^{18})} \Delta_{T^+_{2,3,3,3}}(t^6)$  & $\frac{(1-t^3)^4(1-t^6)(1-t^9)^2}{(1-t^{18})} \Delta_{T_{2,3,3,3}}(t^6)$  \\
(H)=$I$  & $(1-t)^4 \Delta_{T^-_{2,3,5}}(t)$ & $(1-t)^5 \Delta_{T_{2,2,3,5}}(t)$  \\
 (I)  & $(1-t)^3 \Delta_{T^+_{2,3,7}}(t)$ & $(1-t)^4 \Delta_{T_{2,3,4,7}}(t)$   \\
(J)  & $\frac{(1-t^3)^8}{(1-t^6)^2} \Delta_{T^+_{2,2,2,5}}(t^3)$ & $\frac{(1-t^3)^9}{(1-t^6)^2} \Delta_{T_{2,2,2,2,5}}(t^3)$  \\
 (K)  & $\frac{(1-t^3)^6}{(1-t^6)} \Delta_{T^+_{2,4,7}}(t^3)$ & $\frac{(1-t^3)^7}{(1-t^6)} \Delta_{T_{2,2,4,7}}(t^3)$ \\
(L)  & $\frac{(1-t^3)^7}{(1-t^6)} \Delta_{T^+_{2,4,5}}(t^3)$ & $\frac{(1-t^3)^8}{(1-t^6)} \Delta_{T_{2,2,4,5}}(t^3)$  \\
 \hline
\end{tabular}
\end{center}
\caption{Determinants of the matrices $M_0(t)$ and $M(t)$} \label{TabM}
\end{table}

\begin{proof} The character tables of the tetrahedral and icosahedral group are given in \cite{Artin}. The character table of the octahedral group can be found, e.g., in \cite{HH}. From these tables, one can calculate the matrices $B$. The matrices $B$ for the remaining cases are given in \cite{BP}.  For the case (D), only one example is treated. More complete results for the cases $\Delta(3n^2)$ and $\Delta(6n^2)$  can be found in \cite{LNR} and \cite{EL} respectively. From these results, one can derive the corresponding matrices $B$. The proof of the theorem is then obtained by a direct calculation from these matrices using the computer algebra system {\sc Singular} \cite{DGPS}. 
\end{proof}

The results are summarized in Table~\ref{TabM}.

\begin{remark}
Let $G$ be one of the groups $T,O,I$. In this case, the matrix $B$ is symmetric and we have $B^\ast=B$. Therefore 
\[ M(t)=(1-t^3)I-tB+t^2B^\ast = (1-t) ( (1+t+t^2)I -tB).\]
\end{remark}

\section*{Acknowledgements}
This paper grew out of a common research project with A.~Takahashi. The author is grateful to him for useful discussions. He also would like to thank the anonymous referee for useful comments which helped to improve the paper.


\end{document}